\documentclass[reqno]{amsart}
\usepackage{amsmath}
\allowdisplaybreaks[3]
\numberwithin{equation}{section}
\numberwithin{equation}{section}

\def\proof{\indent{\em Proof.\quad}}
\def\endproof{\hfill\hbox{$\sqcup$}\llap{\hbox{$\sqcap$}}\medskip}
\newtheorem{thm}{\indent Theorem}[section]
\newtheorem{cor}[thm]{\indent Corollary}
\newtheorem{lem}[thm]{\indent Lemma}
\newtheorem{prop}[thm]{\indent Proposition}

\newtheorem{rmk}{{\indent\bf Remark}}[section]

\newcommand{\mb}{\mbox}
\newcommand{\hs}{\hspace}

\newcommand{\ima}{\sqrt{-1}}

\newcommand{\td}{\tilde}

\newcommand{\edd}{\end{document}}
\newcommand{\be}{\begin{equation}}
\newcommand{\ee}{\end{equation}}

\newcommand{\lagl}{\langle}
\newcommand{\ragl}{\rangle}
\newcommand{\lmx}{\left(\begin{matrix}}
\newcommand{\rmx}{\end{matrix}\right)}
\newcommand{\ldt}{\left|\begin{matrix}}
\newcommand{\rdt}{\end{matrix}\right|}

\newcommand{\dv}{{\rm div\,}}

\newcommand{\bbr}{{\mathbb R}}
\newcommand{\const}{{\rm const}}

\newcommand{\bbc}{{\mathbb C}}

\newcommand{\ba}{\begin{array}}
\newcommand{\ea}{\end{array}}
\newcommand{\nnm}{\nonumber}
\newcommand{\beal}{\begin{align}}
\newcommand{\eal}{\end{align}}
\newcommand{\bea}{\begin{eqnarray}}
\newcommand{\eea}{\end{eqnarray}}

\newcommand{\spn}{{\rm Span\,}}

\begin{document}

\title[A rigidity theorem of Lagrangian $\xi$-submanifolds in $\mathbb{C}^{2}$]{A rigidity theorem of $\xi$-submanifolds in $\mathbb{C}^{2}$} 

\author[X. X. Li]{Xingxiao Li$^*$} 

\author[X. F. Chang]{Xiufen Chang} 

\dedicatory{}

%
\thanks{Research supported by
National Natural Science Foundation of China (No. 11171091, 11371018).}
\address{
School of Mathematics and Information Sciences
\endgraf Henan Normal University \endgraf Xinxiang 453007, Henan
\endgraf P.R. China}
\email{xxl$@$henannu.edu.cn}

\address{
School of Mathematics and Information Sciences
\endgraf Henan Normal University \endgraf Xinxiang 453007, Henan
\endgraf P.R. China} %
\email{changxfff$@$163.com}



\begin{abstract}
In this paper, we first introduce the concept of $\xi $-submanifold which is a natural generalization of self-shrinkers for the mean curvature flow and also an extension of $\lambda$-hypersurfaces to the higher codimension. Then, as the main results, we prove a rigidity theorem for Lagrangian $\xi $-submanifold in the complex $2$-plane $\bbc^2$.
\end{abstract}

\maketitle

\section{Introduction} 

Let $x: M^n\to\mathbb{R}^{n+p}$ be an $n$-dimensional submanifold in the $(n+p)$-dimensional Euclidean space $\mathbb{R}^{n+p}$. Then $x$ is called a {\em self-shrinker} (to the mean curvature flow) in $\mathbb{R}^{n+p}$ if its mean curvature vector field $H$ satisfies
\be\label{eq1.1}
H+x^{\bot}=0,
\ee
where $x^\bot$ is the orthogonal projection of the position vector $x$ to the normal space $T^\bot M^n$ of $x$.

It is well known that the self-shrinker plays an important role in the study of the mean curvature flow. Not only self-shrinkers correspond to self-shrinking solutions to the mean curvature flow, but also they describe all possible Type I blow ups at a given singularity of the flow. Up to now, there have been a plenty of research papers on self-shrinkers among which are many that provide various results of classification or rigidity theorems. In particular, there are also interesting results about the Lagrangian self-shrinkers in the complex Euclidean $n$-space $\bbc^n$. For example, in \cite{a2}, Anciaux gives new examples of self-shrinking and self-expanding Lagrangian solutions to the mean curvature flow. In \cite{cl1}, the authors classify all Hamiltonian stationary Lagrangian surfaces in the complex plane $\bbc^2$, which are self-similar solutions of the mean curvature flow and, in \cite{cl2}, several rigidity results for Lagrangian mean curvature flow are obtained. As we know, a canonical example of the compact Lagrangian self-shrinker in $\mathbb{C}^2$ is the Clifford torus $\mathbb{S}^{1}(1)\times\mathbb{S}^{1}(1)$.

Recently in \cite{l-w}, Li and Wang prove a rigidity theorem which improves a previous theorem by Castro and Lerma in \cite{cl2}.

\begin{thm}[\cite{l-w}, cf. \cite{cl2}]\label{thm1.1} Let $x:M^2\to \bbr^2$ be a compact orientable Lagrangian self-shrinker with $h$ its second fundamental form. If $|h|^{2}\leq 2$, then $|h|^{2}=2$ and $x(M^{2})$ is the Clifford torus $\mathbb{S}^{1}(1)\times\mathbb{S}^{1}(1)$, up to a holomorphic isometry on $\bbc^2$.
\end{thm}

\begin{rmk}\rm Castro and Lerma also proved Theorem \ref{thm1.1} in \cite{cl2} under the additional condition that the Gauss curvature $K$ of $M^{2}$ is either non-negative or non-positive.\end{rmk}

To make an extension of hypersurface self-shrinkers, Cheng and Wei recently introduce in \cite{C-W2} the definition of $\lambda$-hypersurface of weighted
volume-preserving mean curvature flow in Euclidean space, and classify complete $\lambda$-hypersurfaces
with polynomial area growth and $H-\lambda\geq 0$, which are generalizations of the results
due to Huisken \cite{h} and Colding-Minicozzi \cite{c-m}. According to \cite{C-W2}, a hypersurface $x: M^{n}\to\mathbb{R}^{n+1}$ is called a $\lambda$-hypersurface if its mean curvature $H_0$ satisfies
\be\label{eq1.2}
H_{0}+\langle x, N\rangle=\lambda
\ee
for some constant $\lambda$, where $N$ is the unit normal vector of $x$. Some rigidity or classification results for $\lambda$-hypersurfaces are obtained, for example, in \cite{C-O-W}, \cite{C-W3} and \cite{G2}; For the rigidity theorems for space-like $\lambda$-hypersurfaces see \cite{lc}.

As a natural generalization of both self-shrinkers and $\lambda$-hypersurfaces, we introduce the concept of $\xi$-submanifolds. Precisely, an immersed submanifold $x: M^{n}\to\mathbb{R}^{n+p}$ is called a $\xi$-submanifold if there is a parallel normal vector field $\xi$ such that the mean curvature vector field $H$ satisfies
\be\label{eq1.3}
H+x^{\bot}=\xi.
\ee

Obviously, the Clifford tori $\mathbb{S}^{1}(a)\times\mathbb{S}^{1}(b)$ with positive numbers $a$ and $b$ are examples of Lagrangian $\xi$-submanifold in $\bbc^2$. Similar examples in higher dimensions can be listed as those in \cite{c-l7} for self-shrinkers. In this paper, we focus on the rigidity of compact Lagrangian $\xi$-submanifolds in $\mathbb{C}^2$. Our main theorem is as follows:

\begin{thm}\label{main} Let $x: M^{2}\to\mathbb{C}^{2}$ be a compact orientable Lagrangian $\xi$-submanifold with the second fundamental form $h$ and mean curvature vector $H$. Assume that
$$|h|^{2}+|H-\xi|^{2}\leq |\xi|^{2}+4$$
and $\langle H, \xi\rangle$ is constant. If any one of the following four conditions holds:
\be\label{1.4}
(1)\ |h|^{2}\geq 2,\quad (2)\ |H|^{2}\geq 2,\quad (3)\ |h|^{2}\geq \langle H, H-\xi\rangle,\quad (4)\ \langle H, \xi\rangle\geq 0,
\ee
then
$$|h|^{2}+|H-\xi|^{2}\equiv |\xi|^{2}+4$$
and, up to a holomorphic isometry on $\mathbb{C}^{2}$, $x(M^{2})=\mathbb{S}^{1}(a)\times\mathbb{S}^{1}(b)$ is a standard torus, where $a$ and $b$ are positive numbers satisfying $a^{2}+b^{2}\geq 2a^{2}b^{2}$.
\end{thm}

{\rmk\rm We believe that the last condition \eqref{1.4} in Theorem \ref{main} can be removed. On the other hand, the condition that $\langle H,\xi\rangle$ is constant may also be removed.}

\begin{cor}\label{cor1}
Let $x: M^{2}\to\mathbb{C}^{2}$ be a compact orientable Lagrangian self-shrinker. If
$$|h|^{2}+|H|^{2}\leq 4,$$ then $|h|^{2}+|H|^{2}\equiv 4$ and $x(M^{2})=\mathbb{S}^{1}(1)\times\mathbb{S}^{1}(1)$ up to a holomorphic isometry on $\mathbb{C}^{2}$.
\end{cor}

Clearly, Corollary \ref{cor1} can be viewed as a different new version of Theorem \ref{thm1.1}.

\begin{rmk}\rm Cheng and Wei have introduced in \cite{C-W2} a weighted area functional $\mathcal A$ and derived a related variation formula. Beside the relation between $\lambda$-hypersurfaces and the weighted volume preserving mean curvature flow, they also prove that $\lambda$-hypersurfaces are the critical points of the weighted area functional. Based on this, we believe that similar conclusions will be valid for the $\xi$-submanifolds defined above. Furthermore,
We reasonably believe that, if self-shrinkers and $\lambda$-hypersurfaces take the places of minimal submanifolds and constant mean curvature hypersurfaces, respectively, then $\xi$-submanifolds must take the place of submanifolds of parallel mean curvature vector.\end{rmk}

\section{Lagrangian submanifolds in $\bbc^n$ and their Maslov class}

Let $\mathbb{C}^{n}$ be the complex Euclidean $n$-space with the canonical complex structure $J$. Through out this paper, $x:M^{n}\to \bbc^n$ always denotes an $n$-dimensional Lagrangian submanifold, and $\nabla$, $D$, $\nabla^{\bot}$ denote, respectively, the Levi-Civita connections on $M^{n}$, $\mathbb{C}^{n}$, and the normal connection on the normal boundle $T^\bot M^n$. The formulas of Gauss and Weingarten are given by
$$D_{X}Y=\nabla_{X}Y+h(X, Y),\quad D_{X}\eta=-A_{\eta}X+\nabla^{\bot}_{X}\eta,$$
where $X$, $Y$ are tangent vector fields on $M^n$ and $\eta$ is a normal vector field of $x$.
The Lagrangian condition implies that
$$\nabla^{\bot}_{X}JY=J\nabla_{X}Y, \quad A_{JX}Y=-Jh(X,Y)=A_{JY}X,$$
where $h$ and $A$ are the second fundamental form and the shape operator of $x$, respectively. In particular, $\langle h(X, Y), JZ \rangle$ is totally symmetric as a $3$-form, namely
\be\label{2.1}
\lagl h(X, Y), JZ\ragl= \lagl h(X, Z), JY\ragl= \lagl h(Y, Z), JX\ragl.
\ee

From now on, we agree with the following convention on the ranges of indices:
$$
1\leq i,j,\cdots\leq n,\quad n+1\leq \alpha,\beta,\cdots\leq 2n,\quad 1\leq A,B,\cdots\leq 2n,\quad i^*=n+i.
$$
For a Lagrangian submanifold $x:M^{n}\to\mathbb{C}^{n}$, there are orthonormal frame fields of the form $\{e_i,e_{i^*}\}$ for $\bbc^n$ along $x$, where $e_i\in TM^n$ and $e_{i^*}=Je_i$. Such a frame
is called an {\em adapted Lagrangian frame field} in the literature. The dual frame field is always denoted by $\{\theta_i,\theta_{i^*}\}$, where $\theta_{i^*}=-J\theta_i$. Write
$$h=\sum h^{k^*}_{ij}\theta_i\theta_je_{k^*},\quad \mb{where } h^{k^*}_{ij}=\lagl h(e_i,e_j),e_{k^*}\ragl,$$
or equivalently,
$$h(e_{i}, e_{j})=\sum_{k}h_{ij}^{k^{*}}e_{k^{*}},\quad\mb{ for all }e_i,e_j.$$
Then \eqref{2.1} is equivalent to
\begin{equation}\label{eq2.2}
h_{ij}^{k^{*}}=h_{kj}^{i^{*}}=h_{ik}^{j^{*}}, \quad 1\leq i, j, k\leq n.
\end{equation}
If $\theta_{ij}$ and $\theta_{i^*j^*}$ denote the connection forms of $\nabla$ and $\nabla^\bot$, respectively, then the components $h_{ij,l}^{k^{*}}$, $h_{ij,lp}^{k^{*}}$ of the covariant derivatives of $h$ are given respectively by
\begin{align}
&\sum_{l}h_{ij,l}^{k^{*}}\theta_{l}=dh_{ij}^{k^{*}}+\sum_{l}h_{lj}^{k^{*}}\theta_{li} +\sum_{l}h_{il}^{k^{*}}\theta_{lj}+\sum_{m}h_{ij}^{m^{*}}\theta_{m^{*}k^{*}};\\
&\sum_{p}h_{ij,lp}^{k^{*}}\theta_{p}=dh_{ij,l}^{k^{*}}+\sum_{p}h_{pj,l}^{k^{*}}\theta_{pi} +\sum_{p}h_{ip,l}^{k^{*}}\theta_{pj}+\sum_{p}h_{ij,p}^{k^{*}}\theta_{pl} +\sum_{p}h_{ij,l}^{p^{*}}\theta_{p^{*}k^{*}}.
\end{align}
Moreover, the equations of motion are as follows:
\begin{align}
&dx=\sum_{i}\theta_{i}e_{i},\quad
de_{i}=\sum_{j}\theta_{ij}e_{j}+\sum_{k,j}h_{ij}^{k^{*}}\theta_{j}e_{k^{*}},\\
&de_{k^{*}}=-\sum_{i,j}h_{ij}^{k^{*}}\theta_{j}e_{i}+\sum_{l}\theta_{k^{*}l^{*}}e_{l^{*}}.
\end{align}

Let $R_{ijkl}$ and $R_{i^{*}j^{*}kl}$ denote the components of curvature operators of $\nabla$ and $\nabla^{\bot}$, respectively. Then the equations of Gauss, Codazzi and Ricci are as follows:
\begin{align}
&R_{mijk}=\sum_{l}(h_{mk}^{l^{*}}h_{ij}^{l^{*}}-h_{mj}^{l^{*}}h_{ik}^{l^{*}}), \quad 1\leq m,i, j, k\leq n,\\
&\hs{2cm}h_{ij,l}^{k^{*}}=h_{il,j}^{k^{*}}, \quad 1\leq i, j, k, l\leq n,\label{2.8}\\
&R_{i^{*}j^{*}kl}=\sum_{m}(h_{ml}^{i^{*}}h_{mk}^{j^{*}}-h_{mk}^{i^{*}}h_{ml}^{j^{*}}), \quad 1\leq i, j, k, l\leq n.
\end{align}
The scalar curvature of $\nabla$ is
\begin{equation}
R=|H|^{2}-|h|^{2}\ \mb{ with }|H|^{2}=\sum_{k}\left(\sum_{i}h^{k^*}_{ii}\right)^2,\quad |h|^{2}=\sum_{i,j,k}(h_{ij}^{k^{*}})^{2},
\end{equation}
since the mean curvature vector field
$$H=\sum_{k}H^{k^{*}}e_{k^{*}}=\sum_{i,k} h_{ii}^{k^{*}}e_{k^{*}}.$$

Combining \eqref{eq2.2} and \eqref{2.8}, we know that $h_{ij,l}^{k^{*}}$ is totally symmetric, namely
\begin{equation}\label{eq2.11}
h_{ij,l}^{k^{*}}=h_{jl,k}^{i^{*}}=h_{lk,i}^{j^{*}}=h_{ki,j}^{l^{*}},\quad  1\leq i,j,k,l\leq n,
\end{equation}
and the Ricci identities are as follows:
\begin{equation}
h_{ij,lp}^{k^{*}}-h_{ij,pl}^{k^{*}}=\sum_{m}h_{mj}^{k^{*}}R_{imlp}
+\sum_{m}h_{im}^{k^{*}}R_{jmlp}+\sum_{m}h_{ij}^{m^{*}}R_{k^{*}m^{*}lp}.
\end{equation}

Note that, with respect to the adapted Lagrangian frame $\{e_i,e_{i^*}\}$, the connection forms $\theta_{i^*j^*}=\theta_{ij}$. It follows that
\begin{equation}
R_{m^{*}i^{*}jk}=R_{mijk},\quad\forall m,i,j,k.
\end{equation}

Furthermore, the first and second derivatives $H^{k^*}_{,i}$, $H^{k^*}_{,ij}$ of the mean curvature vector field $H$ are given as
\be
H^{k^*}_{,i}=\sum_jh^{k^*}_{jj,i},\quad H^{k^*}_{,ij}=\sum_lh^{k^*}_{ll,ij}.
\ee

For any smooth function $f$ on $M^{n}$, the covariant derivatives $f_{,i}$, $f_{,ij}$ of $f$, the Laplacian of $f$ are respectively defined as follows:
\begin{equation}
df=\sum_{i}f_{,i}\theta_{i}, \quad \sum_{j}f_{,ij}\theta_{j}=df_{,i}-\sum_{j}f_{,j}\theta_{ij},\quad \triangle f=\sum_{i}f_{,ii}.
\end{equation}

Finally, we also need to introduce the Lagrangian angles, Maslov form and Maslov class of a Lagrangian submanifold in $\bbc^n$ which we shall make use of later.

Let $(z^1,\cdots,z^n)$ be the standard complex coordinates on $\bbc^n$. Then $\Omega=dz^1\wedge\cdots\wedge dz^n$ is a globally defined {\em holomorphic volume form} which is clearly parallel. For a Lagrangian submanifold $x:M^n\to \bbc^n$, {\em the Lagrangian angle} of $x$ is by definition a multi-valued function $\beta:M^n\to \mathbb{R}/2\pi\mathbb{Z}$ given by
$$\Omega_M:=x^*\Omega=e^{\ima\beta}dV_M.$$

As one knows, although the Lagrangian angle $\beta$ can not be determined globally in general, its gradient $\nabla \beta$ is clearly a well-defined vector field on $M^n$, or the same, $\alpha:=d\beta$ is a globally defined $1$-form which is called the {\em Maslov form} of $x$. Clearly, $\alpha$ is closed and thus represents a cohomology class $[\alpha]\in H^1(M^n)$ called the {\em Maslov class}.

In \cite{mo}, the author proved an important formula by which the mean curvature and the Lagrangian angle of a Lagrangian submanifold are linked to each other; A. Arsie has extended this result in \cite{ar} to Lagrangian submanifolds in a general Calabi-Yau manifold.

\begin{thm}[\cite{mo}]\label{thm2.1} Let $x:M^n\to \bbc^n$ be a Lagrangian submanifold and $J$ be the canonical complex structure of $\bbc^n$. Then the mean curvature vector $H$ and the Lagrangian angle $\beta$ meet the following formula:
\be\label{2.16} x_*(\nabla\beta)=-JH.\ee
\end{thm}

\begin{cor}[\cite{s}, \cite{cl2}]\label{cor2.1} Let $x:M^n\to \bbc^n$ be a compact and orientable Lagrangian self-shrinkers. Then the Maslov class $[\alpha]$ can not be trivial. In particular, there does not exist any Lagrangian self-shrinker in $\bbc^n$ with the topology of a sphere.\end{cor}

\begin{rmk}\label{rmk2.1}\rm For our use in this paper, it is necessary to show that Corollary \ref{cor2.1} is still true if we replace the self-shrinker by a $\xi$-submanifold. Precisely, we need\end{rmk}

\begin{prop}\label{prop2.1}
Let $x: M^{n}\to\mathbb{C}^{n}$ be a Lagrangian $\xi$-submanifold. If $M$ is compact and orientable, then $[\alpha]\neq 0$; Consequently, there does not exist any Lagrangian $\xi$-submanifold in $\bbc^n$ with the topology of a sphere.
\end{prop}

\proof By the definition of a $\xi$-submanifold, we have $x=x^{\top}+\xi-H$. By Gauss and Weingarten formulas it follows that, for any $v\in TM^n$,
\begin{align*}
A_{H}v=&-D_{v}H+\nabla^{\bot}_{v}H
=-D_{v}(\xi-x^{\bot})+\nabla^{\bot}_{v}H\\
=&D_{v}x^{\bot}-D_{v}\xi+\nabla^{\bot}_{v}H
=D_{v}x-D_{v}x^{\top}-D_{v}\xi+\nabla^{\bot}_{v}H\\
=&v-\nabla_{v}x^{\top}+A_\xi(v)-h(v,x^{\top})+\nabla^{\bot}_{v}H,
\end{align*}
where $A_{H}$ and $A_\xi$ are Weingarten transformations with respect to $H$ and $\xi$, respectively. Thus
$$A_{H}v=v-\nabla_{v}x^{\top}+A_\xi(v),\quad \nabla_{v}^{\bot}H=h(v,x^{\top}).$$
So that
\begin{align}\label{2.17}
\dv JH=&\sum_{i}\langle\nabla_{e_{i}}JH,e_{i}\rangle
=\sum_{i}\langle J\nabla_{e_{i}}JH, Je_{i} \rangle
=\sum_{i}\langle-\nabla_{e_{i}}^{\bot}H, Je_{i}\rangle\nnm\\
=&\sum_{i}\langle-h(e_{i}, x^{\top}), Je_{i}\rangle
=\sum_{i}-\langle h(e_{i},e_{i}), Jx^{\top}\rangle\nnm\\
=&\sum_{i}\langle Jh(e_{i},e_{i}), x^{\top}\rangle=\langle JH, x^{\top}\rangle,
\end{align}
where $\dv$ is the divergence operator. By \eqref{2.16} and \eqref{2.17} we obtain
\be\label{2.18}
\Delta\beta=\lagl\nabla\beta,x^\top\ragl=\frac{1}{2}\langle\nabla\beta, \nabla|x|^{2}\rangle.
\ee

If $[\alpha]= 0$, then there exists a globally defined Lagrangian angle $\beta$ such that $\alpha=-d\beta$, implying \eqref{2.18} holds globally on $M^n$. Then the compactness assumption and the maximum principle for a second linear elliptic partial equation (see \cite{gt}, for example) assure that $\beta$ must be constant. Hence $H=x_*(J\nabla\beta)\equiv 0$, contradicting to the fact that there are no compact minimal submanifolds in Euclidean space. This contradiction proves that $[\alpha]\neq 0$.

Since the first homology of a sphere $S^n$ vanishes for $n>1$, there can not be any Lagrangian $\xi$-submanifolds with the topology of a sphere.\endproof

\section{Proof of the main theorem}

Let $x: M^n\to \mathbb{C}^n$ be a Lagrangian $\xi$-submanifold without boundary. Then, with respect to a orthonormal frame field $\{e_i\}$, the defining equation \eqref{eq1.3} is equivalent to
\begin{equation}\label{eq3.1}
H^{k^{*}}=-\langle x, e_{k^{*}}\rangle+\xi^{k^{*}}, \quad 1 \leq k \leq n.
\end{equation}
where $\xi=\sum \xi^{k^*}e_{k^*}$ is a given parallel normal vector field. From now on, we always assume that $n=2$ if no other specification is given.

We start with a well-known operator $\mathcal{L}$ acting on smooth functions defined by
\be\label{3.4}\mathcal{L}=\triangle -\langle x, \nabla\cdot\rangle =e^{\frac{|x|^{2}}2}\dv(e^{-\frac{|x|^{2}}{2}}\nabla\cdot),\ee
which was first introduced by Colding and Minicozzi in \cite{c-m} to the study of self-shrinkers. Since then, the operator $\mathcal{L}$ has been one of the most effect tools adapted by many authors. In particular, the following is a fundamental lemma related to $\mathcal{L}$:

\begin{lem} [\cite{lw2}]\label{lem3.1} Let $x: M^n\to \mathbb{R}^{n+p}$ be a complete immersed submanifold. If $u$ and $v$ are $C^2$-smooth functions with
$$\int_M(|u\nabla v|+|\nabla u\nabla v|+|u\mathcal{L} v|)e^{-\frac{|x|^2}{2}}dV_M< \infty,$$
then it holds that
$$\int_M u\mathcal{L}v e^{-\frac{|x|^2}{2}}dV_M=-\int_M\lagl \nabla u, \nabla v\ragl e^{-\frac{|x|^2}{2}}dV_M.$$\end{lem}

Now, to make the whole argument more readable, we divide our proof into the following lemmas and propositions:

\begin{lem}[cf.\cite{l-w}]\label{lem3-2}
Let $x: M^{2}\to \mathbb{C}^{2}$ be a Lagrangian $\xi$-submanifold. Then
\begin{align}
&H_{,i}^{k^{*}}=\sum_{j}h_{ij}^{k^{*}}\langle x, e_{j}\rangle, \quad 1\leq i,k\leq 2,\label{eq3.2}
\\
&H_{,ij}^{k^{*}}=\sum_{m}h_{im,j}^{k^{*}}\langle x, e_{m}\rangle+h_{ij}^{k^{*}}-\sum_{m,p}(H-\xi)^{p^{*}}h_{im}^{k^{*}}h_{mj}^{p^{*}},\quad 1\leq i,j,k\leq 2.
\label{eq3.3}\end{align}
\end{lem}

\begin{lem}\label{lem3.2} It holds that
\begin{align}
\frac{1}{2}\mathcal{L}(|h|^{2}+|H-\xi|^{2})=&|\nabla h|^{2}+|\nabla^{\bot}H|^{2}+|h|^{2}\nnm\\ &-\frac{1}{2}(|h|^{2}-|H|^{2})(3|h|^{2}-2|H|^{2}+\langle H, H-\xi\rangle)\nnm\\
&+\langle H,H-\xi\rangle -\sum_{i,j,k,l}h_{ij}^{k^{*}}h_{ij}^{l^{*}}(H-\xi)^{k^{*}}(H-\xi)^{l^{*}}\nnm\\
&-\sum_{i,j,k,l}h_{ij}^{k^{*}}h_{ij}^{l^{*}}H^{k^{*}}(H-\xi)^{l^{*}}.\label{3.5}
\end{align}
\end{lem}

\proof By a direct computation using Lemma \ref{lem3-2} we find (cf. \cite{l-w})
\begin{align}
\frac{1}{2}\mathcal{L}|h|^{2}=&|\nabla h|^{2}+|h|^{2}-\frac{3}{2}|h|^{4} +\frac{5}{2}|H|^{2}|h|^{2}-|h|^{4}\nnm\\
&+\frac{1}{2}\langle H, H-\xi\rangle(|H|^{2}-|h|^{2})-
\sum_{i,j,k,l}H^{k^{*}}h_{ij}^{k^{*}}h_{ij}^{l^{*}}(H-\xi)^{l^{*}};\\
\frac{1}{2}\mathcal{L}(|H-\xi|^{2})=&\frac{1}{2}\triangle(|H-\xi|^{2}) -\frac{1}{2}\langle x, \nabla|H-\xi|^{2}\rangle\nnm\\
=&\sum_{i,k}(H-\xi)^{k^{*}}H_{,ii}^{k^{*}}+|\nabla^{\bot}H|^{2}
-\sum_{i.k}(H-\xi)^{k^{*}}H_{,i}^{k^{*}}\langle x, e_{i}\rangle\nnm\\
=&\langle H- \xi, H\rangle+|\nabla^{\bot}H|^{2}
-\sum_{i,j,k,l}(H-\xi)^{k^{*}}h_{ij}^{k^{*}}h_{ij}^{l^{*}}(H-\xi)^{l^{*}}.
\end{align}
By taking the sum we obtain \eqref{3.5}.\endproof

\begin{lem}\label{lem3.3} It holds that
\begin{align}
\frac{1}{2}\triangle(|x^{\top}|^{2})=&\sum_{i,j,k}h_{ij}^{k^{*}}\langle x,e_{i}\rangle \langle x, e_{j}\rangle(\xi-H)^{k^{*}}-\sum_{i,j,k,l}h_{il}^{k^{*}}h_{lj}^{k^{*}}
\langle x, e_{i}\rangle\langle x, e_{j}\rangle\nnm\\
&+2-2\langle H,H-\xi\rangle +\sum_{i,j,k,l}h_{ij}^{k^{*}}h_{ij}^{l^{*}}(H-\xi)^{k^{*}}(H-\xi)^{l^{*}}.
\label{3.8}
\end{align}
\end{lem}

\proof We find
\begin{align*}
\frac{1}{2}\triangle(|x^{\top}|^{2})=&\frac{1}{2}\sum_{i,j}\langle x, e_{j}\rangle^{2}_{,ii}=\sum_{i,j}(\langle x, e_{j}\rangle\langle x, e_{j}\rangle_{i})_{,i}\\
=&\sum_{i,j,k}(\langle x, e_{j}\rangle\langle x_{i}, e_{j}\rangle+\langle x, e_{j}\rangle\langle x, h_{ji}^{k^{*}}e_{k^{*}}\rangle)_{,i}\\
=&2+2\sum_{i,k}h_{ii}^{k^{*}}\langle x, e_{k^{*}}\rangle+\sum_{j,k}H_{,j}^{k^{*}}\langle x, e_{j}\rangle\langle x, e_{k^{*}}\rangle\\
&+\sum_{i,j,k,l}h_{ij}^{k^{*}}h_{ij}^{l^{*}}\langle x, e_{l^{*}}\rangle\langle x, e_{k^{*}}\rangle-\sum_{i,j,k,l}h_{ij}^{k^{*}}h_{il}^{k^{*}}\langle x, e_{j}\rangle\langle x, e_{l}\rangle\\
=&2-2\langle H,H-\xi\rangle+\sum_{i,j,k}h_{ij}^{k^{*}}\langle x, e_{i}\rangle\langle x, e_{j}\rangle(\xi-H)^{k^{*}}\\
&+\sum_{i,j,k,l}h_{ij}^{k^{*}}h_{ij}^{l^{*}}(H-\xi)^{l^{*}}(H-\xi)^{k^{*}}
-\sum_{i,j,k,l}h_{il}^{k^{*}}h_{lj}^{k^{*}}\langle x, e_{i}\rangle\langle x, e_{j}\rangle,
\end{align*}
and the lemma is proved.\endproof

\begin{lem}\label{lem3.4} It holds that
\begin{align}
&\triangle(\langle H, \xi\rangle)=\sum_{i,j,k}h_{ij}^{k^{*}}\langle x, e_{i}\rangle\langle x, e_{j}\rangle\xi^{k^{*}}+\langle H, \xi\rangle-\sum_{i,j,k,l}h_{ij}^{k^{*}}h_{ij}^{l^{*}}\xi^{k^{*}}(H-\xi)^{l^{*}},\\
&\mathcal{L}(\langle H, \xi\rangle)=\langle H, \xi\rangle-\sum_{i,j,k,l}h_{ij}^{k^{*}}h_{ij}^{l^{*}}\xi^{k^{*}}(H-\xi)^{l^{*}}.\label{3.10}
\end{align}
\end{lem}

\proof By \eqref{eq3.2} and \eqref{eq3.3},
\begin{align*}
\triangle(\langle H, \xi\rangle)=&\sum_{i,k}(H^{k^{*}}\xi^{k^{*}})_{,ii}=\sum_{i,k}H^{k^{*}}_{,ii}\xi^{k^{*}}\\
=&\sum_{i,k,l,m}(h_{im,i}^{k^{*}}\langle x, e_{m}\rangle+h_{ii}^{k^{*}}-(H-\xi)^{^{l^{*}}}h_{im}^{k^{*}}h_{mi}^{l^{*}})\xi^{k^{*}}\\
=&\sum_{i,k}H_{,i}^{k^{*}}\langle x, e_{i}\rangle\xi^{k^{*}}+\langle H, \xi\rangle-
\sum_{i,j,k,l}h_{ij}^{k^{*}}h_{ij}^{l^{*}}\xi^{k^{*}}(H-\xi)^{l^{*}}\\
=&\sum_{i,j,k}h_{ij}^{k^{*}}\langle x, e_{i}\rangle\langle x, e_{j}\rangle\xi^{k^{*}}+\langle H, \xi\rangle-\sum_{i,j,k,l}h_{ij}^{k^{*}}h_{ij}^{l^{*}}\xi^{k^{*}}(H-\xi)^{l^{*}};\\
\langle x, \nabla\langle H, \xi\rangle\rangle=&\sum_{i}\langle H, \xi\rangle_{,i}\langle x, e_{i}\rangle=\sum_{i,j,k}h_{ij}^{k^{*}}\langle x, e_{i}\rangle\langle x, e_{j}\rangle\xi^{k^{*}}.
\end{align*}
Thus, by adding them up, we get \eqref{3.10}. \endproof

\begin{lem}[ cf. \cite{c-l7}, \cite{c-m}; also \cite{l-w}]\label{lem3.5} It holds that
\begin{align*}&\frac{1}{2}\triangle(|x|^{2})=2-\langle H, H-\xi\rangle,\\
\frac{1}{2}\mathcal{L}(&|x|^{2})=|\xi|^{2}+2-(|x|^{2}+\langle H, \xi\rangle).\end{align*}
\end{lem}

\proof From \eqref{eq3.1}, we find
\begin{align*}
\frac{1}{2}\triangle(|x|^{2})=&2+\langle x, \triangle x\rangle
=2+\sum_{k}H^{k^{*}}\langle x, e_{k^{*}}\rangle
=2-\langle H,H-\xi\rangle,\\
\frac{1}{2}\mathcal{L}(|x|^{2})=&\frac{1}{2}\triangle(|x|^{2})-\frac{1}{2}\langle x, \nabla|x|^{2}\rangle=2-|H|^{2}+\langle H,\xi\rangle-|x^{\top}|^{2}\\
=&2+|\xi|^{2}-(|x|^{2}+\langle H, \xi\rangle).
\end{align*}
\endproof

\begin{prop}\label{prop3.6} Let $M^2$ be orientable and compact. If
\begin{align*}
|h|^{2}+|H-\xi|^{2}\leq |\xi|^{2}+4,
\end{align*}
then
\be\label{3.11}|h|^{2}+|H-\xi|^{2}\equiv |\xi|^{2}+4\ee
and $x(M^{2})$ is a topological torus.
\end{prop}

\proof By Lemma \ref{lem3.5},
\begin{align}
\int_{M}|H-\xi|^{2}dV_{M}=&\int_{M}(|\xi|^{2}+2(|H|^{2}-\langle H, \xi\rangle)-|H|^{2})dV_{M}\nnm\\
=&\int_{M}(|\xi|^{2}+4-|H|^{2})dV_{M}.\label{3.12}
\end{align}
Let $K$ be the Gauss curvature of $M^2$. Then the Gauss equation gives that
$$2K= |H|^{2}-|h|^{2}.$$
Denote by ${\rm gen}(M^2)$ the genus of $M^2$. Then from the Gauss-Bonnet theorem and \eqref{3.12} it follows that
\begin{align*}
8\pi(1-{\rm gen}(M^2))=&2\int_{M}KdV_{M}=\int_{M}(|H|^{2}-|h|^{2})dV_{M}\\
=&\int_{M}(|\xi|^{2}+4-(|h|^{2}+|H-\xi|^{2}))dV_{M}\geq 0,
\end{align*}
implying that ${\rm gen}(M^2)\geq 1$. So $M^2$ is topologically either a $2$-sphere or a torus. But Proposition \ref{prop2.1} excludes the first possibility. So ${\rm gen}(M^2)=1$ and \eqref{3.11} is proved. \endproof

\begin{lem}\label{lem3.8}
Let $p_{0}\in M^2$ be a point where $|x|^{2}$ attains its minimum on $M^2$. If $M^2$ is orientable, compact and
$$|h|^{2}+|H-\xi|^{2}=\const,$$
then
\be\label{3.13} \nabla^{\bot}H(p_{0})=0, \quad(\nabla h)(p_{0})=0.\ee
\end{lem}

\proof Since $(|x|^{2})_{,j}=0,1\leq j\leq 2$ at $p_{0}$, it holds that $\langle x, e_{j}\rangle(p_{0})=0$, $1\leq j\leq 2$. So by \eqref{eq3.2} we have
\be\label{3.14}
H^{k^{*}}_{,i}=0, \ 1\leq i,k \leq 2,\quad
|H-\xi|^{2}_{,i}=\sum_{k}(H-\xi)^{k^{*}}H_{,i}^{k^{*}}
=0, \quad 1\leq i\leq 2 \quad \mb{at\ }\ p_{0}
\ee
where the first set of equalities are exactly $\nabla^\bot H(p_0)=0$, which give
\begin{equation}\label{eq3.9}
h_{11,1}^{1^{*}}+h_{22,1}^{1^{*}}=0,\quad h_{11,2}^{1^{*}}+h_{22,2}^{1^{*}}=0,\quad
h_{11,1}^{2^{*}}+h_{22,1}^{2^{*}}=0,\quad
h_{11,2}^{2^{*}}+h_{22,2}^{2^{*}}=0.
\end{equation}

On the other hand, from
\begin{equation}\label{eq3.5}
|h|^{2}+|H-\xi|^{2}=\const,
\end{equation}
we obtain
\begin{equation}\label{eq3.6}
|h|^{2}_{,k}+|H-\xi|^{2}_{,k}\equiv0,\quad 1\leq k\leq 2,
\end{equation}
which with \eqref{3.14} implies that
$$(|h|^{2})_{,k}=0,\quad 1\leq k\leq 2\quad \mb{at\ }\ p_{0}.$$
Since
$$|h|^{2}=(h_{11}^{1^{*}})^{2}+2(h_{12}^{1^{*}})^{2}+(h_{22}^{1^{*}})^{2}
+(h_{11}^{2^{*}})^{2}+2(h_{12}^{2^{*}})^{2}+(h_{22}^{2^{*}})^{2},$$
we find that
\begin{align}
&h_{11}^{1^{*}}h_{11,1}^{1^{*}}+2h_{12}^{1^{*}}h_{12,1}^{1^{*}}+h_{22}^{1^{*}}h_{22,1}^{1^{*}}
+h_{11}^{2^{*}}h_{11,1}^{2^{*}}+2h_{12}^{2^{*}}h_{12,1}^{2^{*}}
+h_{22}^{2^{*}}h_{22,1}^{2^{*}}=0,\label{eq3.10}\\
&h_{11}^{1^{*}}h_{11,2}^{1^{*}}+2h_{12}^{1^{*}}h_{12,2}^{1^{*}}+h_{22}^{1^{*}}h_{22,2}^{1^{*}}
+h_{11}^{2^{*}}h_{11,2}^{2^{*}}+2h_{12}^{2^{*}}h_{12,2}^{2^{*}}
+h_{22}^{2^{*}}h_{22,2}^{2^{*}}=0\label{eq3.11}
\end{align}
hold at $p_{0}$.
From \eqref{eq2.11} and \eqref{eq3.9} we get
\begin{equation}\label{eq3.12}
h_{22,1}^{1^{*}}=-h_{11,1}^{1^{*}},\quad h_{22,2}^{1^{*}}=-h_{11,2}^{1^{*}},\quad
h_{22,2}^{2^{*}}=h_{11,1}^{1^{*}}\quad \mb{at\ }\ p_{0}.
\end {equation}

Since, by \eqref{eq2.2} and \eqref{eq2.11}, both $h_{ij}^{k^{*}}$ and $h_{ij,l}^{k^{*}}$ are totally symmetric, we obtain by \eqref{eq3.12}, \eqref{eq3.10} and \eqref{eq3.11} that
\begin{align}
&(h_{11}^{1^{*}}-3h_{22}^{1^{*}})h_{11,1}^{1^{*}}-
(h_{22}^{2^{*}}-3h_{11}^{2^{*}})h_{11,2}^{1^{*}}=0,\label{3.21}\\
&(h_{22}^{2^{*}}-3h_{11}^{2^{*}})h_{11,1}^{1^{*}}+
(h_{11}^{1^{*}}-3h_{22}^{1^{*}})h_{11,2}^{1^{*}}=0\quad \mb{at\ }\ p_{0}.\label{3.22}
\end{align}

We claim that
\be\label{3.24}(\nabla h)(p_{0})=0.\ee

Otherwise, we should have $(h_{11,1}^{1^{*}})^{2}+(h_{11,2}^{1^{*}})^{2}\neq 0$ at $p_0$. Then from \eqref{3.21} and \eqref{3.22} it follows that $$(h_{11}^{1^{*}}-3h_{22}^{1^{*}})^{2}+(h_{22}^{2^{*}}-3h_{11}^{2^{*}})^{2}=0\quad \mb{at\ }\ p_{0}.$$
Thus
\be
 |h|^{2}(p_{0})=\frac{4}{3}((h_{11}^{1^{*}})^{2}+(h_{22}^{2^{*}})^{2}),\quad
 |H|^{2}(p_{0})=\frac{16}{9}((h_{11}^{1^{*}})^{2}+(h_{22}^{2^{*}})^{2}).
\ee

Now by the definition of $p_0$ and Lemma \ref{lem3.5},
$$0\leq \frac{1}{2}\triangle {|x|^{2}}(p_{0})
=2-\langle H, H-\xi\rangle(p_{0}).$$

It follows that
\begin{align}
|h|^{2}+|H-\xi|^{2}=&(|h|^{2}+|H-\xi|^{2})(p_{0})\nnm\\
=&\frac{3}{4}|H|^{2}(p_{0})+2\langle H, H-\xi\rangle(p_{0})
-|H|^{2}(p_{0})+|\xi|^{2}\nnm\\
=&|\xi|^{2}+2\langle H, H-\xi\rangle(p_{0})-\frac{1}{4}|H|^{2}(p_{0})\\
\leq&|\xi|^{2}+4.\label{3.25}
\end{align}
Therefore, by Proposition \ref{prop3.6}, $|h|^{2}+|H-\xi|^{2}=|\xi|^{2}+4$. But it is easy to see that the equality in \eqref{3.25} holds if and only if $|H|^{2}(p_{0})=0$ and
$\langle H, H-\xi\rangle(p_{0})=2$, which is of course not possible! This contradiction proves the above claim and completes the proof of Lemma \ref{lem3.8}. \endproof

{\rmk\rm Our main observation here is that, if $p_{0}\in M^{2}$ is a minimum point of $|x|^{2}$ then $$x^{\top}(p_{0})=\sum_{i}\langle x, e_{i}\rangle e_{i}(p_{0})=0,$$
implying
$$\nabla^{\bot}H(p_{0})=\nabla^{\bot}(H-\xi)(p_{0})=0.$$
In particular, $p_{0}$ is also a minimum point of $|x^{\top}|^{2}$.}

\begin{prop}\label{prop3.9}
Let $x: M^{2}\to\mathbb{C}^{2}$ be a compact and orientable Lagrangian $\xi$-submanifold. Suppose that $$|h|^{2}+|H-\xi|^{2}=|\xi|^{2}+4$$
and $\langle H, \xi\rangle$ is constant. If any one of the followings holds,
\be\label{eq3.17}
(1)\ |h|^{2}\geq 2,\quad (2)\ |H|^{2}\geq 2,\quad (3)\ |h|^{2}\geq \langle H, H-\xi\rangle,\quad (4)\ \langle H, \xi\rangle\geq 0,
\ee
then $|x|^{2}$ is a constant.
\end{prop}

\proof As above, let $p_0$ be a minimum point of $|x|^2$. Then, by Lemma \ref{lem3.2} and Lemma \ref{lem3.8}, it holds at $p_{0}$ that
\begin{align}
0=&\frac{1}{2}\mathcal{L}(|h|^{2}+|H-\xi|^{2})\nnm\\
=&|h|^{2}-\frac{1}{2}(|h|^{2}-|H|^{2})(3|h|^{2}-2|H|^{2}
+\langle H,H-\xi\rangle)+\langle H,H-\xi\rangle\nnm\\
&-\sum_{i,j,k,l}h_{ij}^{k^{*}}h_{ij}^{l^{*}}(H-\xi)^{k^{*}}(H-\xi)^{l^{*}}
-\sum_{i,j,k,l}h_{ij}^{k^{*}}h_{ij}^{l^{*}}H^{k^{*}}(H-\xi)^{l^{*}}.
\end{align}

Furthermore, form Lemma \ref{lem3.3} and Lemma \ref{lem3.4} it follows that, at $p_{0}$
\begin{align}
&0\leq \frac{1}{2}\triangle(|x^{\top}|^{2})=2-2\langle H,H-\xi\rangle+
\sum_{i,j,k,l}h_{ij}^{k^{*}}h_{ij}^{l^{*}}(H-\xi)^{k^{*}}(H-\xi)^{l^{*}},\\
&0=\frac{1}{2}\mathcal{L}(\langle H,\xi\rangle)
=\frac{1}{2}(\langle H,\xi\rangle-
\sum_{i,j,k,l}h_{ij}^{k^{*}}h_{ij}^{l^{*}}\xi^{k^{*}}(H-\xi)^{l^{*}}),
\end{align}
implying
\be
-\sum_{i,j,k,l}h_{ij}^{k^{*}}h_{ij}^{l^{*}}(H-\xi)^{k^{*}}(H-\xi)^{l^{*}}
\leq 2-2\langle H,H-\xi\rangle,
\ee
and
\begin{align*}
-\sum_{i,j,k,l}h_{ij}^{k^{*}}h_{ij}^{l^{*}}(H-\xi)^{l^{*}}H^{k^{*}}
=&-\sum_{i,j,k,l}h_{ij}^{k^{*}}h_{ij}^{l^{*}}(H-\xi)^{k^{*}}(H-\xi)^{l^{*}}\\
&-\sum_{i,j,k,l}h_{ij}^{k^{*}}h_{ij}^{l^{*}}\xi^{k^{*}}(H-\xi)^{l^{*}} \\
\leq&2-2\langle H,H-\xi\rangle-\langle H,\xi\rangle \\
=&2-\langle H,H-\xi\rangle-|H|^{2}.
\end{align*}

Consequently, we have at $p_{0}$
\begin{align*}
0=&\frac{1}{2}\mathcal{L}(|h|^{2}+|H-\xi|^{2})\\
\leq&-\frac{1}{2}(|h|^{2}-|H|^{2})(3|h|^{2}-2|H|^{2}
+\langle H,H-\xi\rangle)\\
&+|h|^{2}-|H|^{2}+2(2-\langle H,H-\xi\rangle).
\end{align*}

On the other hand, from
$$|h|^{2}+|H-\xi|^{2}=|\xi|^{2}+4,$$
we know that
$$|h|^{2}-|H|^{2}=2(2-\langle H,H-\xi\rangle)\geq 0\quad \mb{at\ }\ p_{0}.$$
Thus, if one of \eqref{eq3.17} holds, then at $p_0$
\begin{align}
0=&\frac{1}{2}\mathcal{L}(|h|^{2}+|H-\xi|^{2})\nnm\\
\leq&-\frac{1}{2}(|h|^{2}-|H|^{2})(2|h|^{2}-|H|^{2}+4
-\langle H,H-\xi\rangle)+2(|h|^{2}-|H|^{2})\nnm\\
=&-\frac{1}{2}(|h|^{2}-|H|^{2})(2|h|^{2}
-|H|^{2}-\langle H,H-\xi\rangle)\nnm\\
=&-\frac{1}{2}(|h|^{2}-|H|^{2})(|h|^{2}-|H|^{2}+|h|^{2}
-\langle H,H-\xi\rangle)\nnm\\
=&-\frac{1}{2}(|h|^{2}-|H|^{2})(|h|^{2}-|H|^{2}+|h|^{2}
-2+2-\langle H,H-\xi\rangle)\nnm\\
=&-\frac{1}{2}(|h|^{2}-|H|^{2})(2(|h|^{2}-|H|^{2})
+\langle H,\xi\rangle)\nnm\\
=&-\frac{1}{2}(|h|^{2}-|H|^{2})(2(|h|^{2}-|H|^{2})+|H|^{2}
-2+2-\langle H,H-\xi\rangle)\nnm
\leq0\nnm.
\end{align}
Consequently
$$|h|^{2}-|H|^{2}=2-\langle H,H-\xi\rangle=0\quad \mb{at\ }\ p_{0}.$$
It follows that
\begin{align}
|x|^{2}+\langle H,\xi\rangle\geq& |x|^{2}(p_{0})+\langle H,\xi\rangle(p_{0})\nnm\\
=&|H-\xi|^{2}(p_{0})+\langle H,\xi\rangle(p_{0})\nnm\\
=&\langle H,H-\xi\rangle(p_{0})+|\xi|^{2}\nnm\\
=&|\xi|^{2}+2\nnm.
\end{align}

This together with Lemma \ref{lem3.1} (for $u=1$, $v=|x|^2$) and \ref{lem3.5} gives that
$$
0=\int_{M}\frac{1}{2}\mathcal{L}(|x|^{2})e^{-\frac{|x|^{2}}{2}}dV_{M}
=\int_{M}(|\xi|^{2}+2-(|x|^{2}+\langle H,\xi\rangle) )e^{-\frac{|x|^{2}}{2}}dV_{M}\leq 0
$$
implying that $|x|^{2}+\langle H,\xi\rangle=|\xi|^{2}+2$. In particular, $|x|^{2}=\const $. \endproof

\begin{prop}\label{prop3.10}
Let $x: M^{n}\to N^{n}$ be a Lagrangian submanifold in a $K\ddot{a}hler$ manifold $N^{n}$. If both $M^{n}$ and $N^{n}$ are flat, then around each point $p\in M^{n}$, there
exists some orthonormal frame field $\{e_{i},e_{i^*}\}$ with $e_{i^{*}}=Je_{i}$ ($1\leq i\leq n$), such that
$$h_{ij}^{k^{*}}:=\langle h(e_{i},e_{j}),e_{k^{*}}\rangle=\lambda_{i}^{k^{*}}\delta_{ij}, \quad 1\leq i,j,k\leq n.$$
\end{prop}

\proof For $p\in M^{n}$, we pick an orthonormal tangent frame $\{\bar{e}_{i}\}$ and an orthonormal normal frame $\{\bar{e}_{\alpha}\}_{n+1\leq \alpha\leq 2n}$. Define
$$\bar{h}_{ij}^{\alpha}=\langle h(\bar{e}_{i},\bar{e}_{j}),\bar{e}_{\alpha}\rangle.$$

Since both $M^n$ and $N$ are flat, $T^{\bot}M^{n}$ is also flat with respect to the normal connection. By the Ricci equation
$$
0=\langle R^{^{\bot}}\langle e_{i},e_{j}\rangle e_{\alpha},e_{\beta}\rangle=\sum _k(h_{ik}^{\alpha}h_{jk}^{\beta}-h_{ik}^{\beta}h_{jk}^{\alpha}).
$$
Hence we can choose another orthonormal tangent frame $\{e_{i}\}$ such that
$$h_{ij}^{\alpha}:=\langle h(e_{i},e_{j}),e_{\alpha}\rangle=\mu_{i}^{\alpha}\delta_{ij}.$$

Write $e_{k^{*}}=\sum_{\alpha}a_{k^{*}}^{\alpha}e_{\alpha}$. Then
\begin{align}
h_{ij}^{k^{*}}=&\langle h(e_{i},e_{j}),e_{k^{*}}\rangle
=\big\langle h(e_{i},e_{j}),\sum_{\alpha}a_{k^{*}}^{\alpha}e_{\alpha}\big\rangle
=\sum_{\alpha}a_{k^{*}}^{\alpha}\langle h(e_{i},e_{j}),e_{\alpha}\rangle\nnm\\
=&\sum_{\alpha}a_{k^{*}}^{\alpha}h_{ij}^{\alpha}
=\sum_{\alpha}a_{k^{*}}^{\alpha}\mu_{i}^{\alpha}\delta_{ij}
=\lambda_{i}^{k^{*}}\delta_{ij}\ \mb{ with } \lambda_{i}^{k^{*}}:=\sum_{\alpha}a_{k^{*}}^{\alpha}\mu_{i}^{\alpha}.
\end{align}
Thus Proposition \ref{prop3.10} is proved.\endproof

{\em Proof of the main theorem:}

Since, by the Gauss equation
$2K=|H|^{2}-|h|^{2}\equiv 0$, namely, $M^{2}$ is flat. Therefore,
due to Proposition \ref{prop3.10}, we can choose $\{e_{1},e_{2}\}$ such that
\be
h_{12}^{1^{*}}=h_{12}^{2^{*}}=0.
\ee
It follows that
\be
h_{22}^{1^{*}}=h_{11}^{2^{*}}=0.
\ee

On the other hand, since $\nabla h\equiv 0$, we have
$$0=\sum h _{ijl}^{k^{*}}\theta_{l}=dh_{ij}^{k^{*}}-\sum h_{lj}^{k^{*}}\theta_{il}
-\sum h_{il}^{k^{*}}\theta_{jl}+\sum h_{ij}^{p^{*}}\theta_{p^{*}k^{*}}.
$$
It follows that
\begin{enumerate}
\item $i=j=1$, $k=2$, we get $ h_{11}^{1^{*}}\theta_{1^{*}2^{*}}=0$,
\item $i=j=2$, $k=1$, we get $h_{22}^{2^{*}}\theta_{2^{*}1^{*}}=0$,
\item $i=j=k=1$, $0=dh_{11}^{1^{*}}+h_{11}^{1^{*}}\theta_{1^{*}1^{*}}=dh_{11}^{1^{*}}$, we get
  $h_{11}^{1^{*}}=\const$,
\item $i=j=k=2$, $0=dh_{22}^{2^{*}}+h_{22}^{2^{*}}\theta_{2^{*}2^{*}}=dh_{22}^{2^{*}}$, we get
  $h_{22}^{2^{*}}=\const$.
\end{enumerate}

Since $x$ can not be totally geodestic, $(h_{11}^{1^{*}})^{2}+(h_{22}^{2^{*}})^{2}\neq 0$. Thus, we have $\theta_{1^{*}2^{*}}=0$, which with $\nabla J=0$ shows that $\theta_{12}=0$. From this it is easy to see that $h_{11}^{1^{*}}h_{22}^{2^{*}}\neq 0$. Let
$$\td e_{1}=e_{1}\cos\theta-e_{2}\sin\theta,
\quad \td e_{2}=e_{1}\sin\theta+e_{2}\cos\theta
$$
be another frame field such that
$$\td h_{ij}^{k^{*}}:=\langle h(\td e_{i},\td e_{j}),\td e_{k^{*}}\rangle
=\td {\lambda}_{i}^{k^{*}}\delta_{ij}.
$$
Then a direct computation shows that
$$
\sin\theta\cos\theta(h_{11}^{1^{*}}\cos\theta+h_{22}^{2^{*}}\sin\theta)
=\sin\theta\cos\theta(h_{11}^{1^{*}}\sin\theta-h_{22}^{2^{*}}\cos\theta)=0.
$$

Since $(h_{11}^{1^{*}})^{2}+(h_{22}^{2^{*}})^{2}\neq 0$, we have $\sin2\theta=0$, that is
$\theta=0,$ or $\frac{\pi}{2}$ or $\pi$. Clealy, by choosing $\theta=\frac{\pi}{2}$, we can change the sign of $h_{11}^{1^{*}}h_{22}^{2^{*}}$; while by choosing $\theta=\pi$, we can change the sign of both $h_{11}^{1^{*}}$ and $h_{22}^{2^{*}}$. Thus we can always assume that $h_{11}^{1^{*}}>0$ and $h_{22}^{2^{*}}> 0$. It then follows that $\{e_{1},e_{2}\} $ can be uniquely determined and, in particular, is globally defined.

Define
\be\label{3.35} V_{1}=\spn_\bbr\{e_{1},e_{1^*}\}=\spn_{\mathbb{C}}\{e_{1}\}, \quad V_{2}=\spn_\bbr\{e_2,e_{2^*}\}=\spn_{\mathbb{C}}\{e_{2}\}.
\ee
Since
\begin{align*}
&de_{1}=\nabla e_{1}+\sum h_{1j}^{k^{*}}\theta_{j}e_{k^{*}}=h_{11}^{1^{*}}\theta_{1}e_{1^{*}}\in V_{1},
\\
&de_{1^{*}}=Jde_{1}=h_{11}^{1^{*}}\theta_{1}Je_{1^{*}}=-h_{11}^{1^{*}}\theta_{1}e_{1}\in V_{1},
\end{align*}
we know that $V_{1}$ is a $1$-dimensional constant complex subspace of $\mathbb{C}^{2}$.

Similarly, $V_{2}$ is also a $1$-dimensional constant complex subspace of $\mathbb{C}^{2}$. Furthermore, $V_1$ and $V_2$ are clearly orthogonal. So, up to a holomorphic isometry on $\mathbb{C}^{2}$, we can assume that $V_{1}=\mathbb{C}^{1}$, $V_{2}=\mathbb{C}^{1}$ so that $\mathbb{C}^{2}=V_{1}\times V_{2}$. Write
$$x=(x^1,x^2)\in V_1\times V_2\equiv \bbc^2.$$
Then
$$0=e_{i}(|x|^{2})=e_{i}(|x^1|^{2})+e_{i}(|x^2|^{2}),\quad i=1,2,$$
which with the definitions \eqref{3.35} of $V_1$ and $V_2$ shows that
$$e_{1}(|x^{1}|^{2})=e_{2}(|x^{1}|^{2})=0,\quad e_{1}(|x^{2}|^{2})=e_{2}(|x^{2}|^{2})=0,$$
that is, $$|x^{1}|^{2}=\const,\quad |x^{2}|^{2}=const.$$
It is easily seen that both $|x^{1}|^{2}$ and $|x^{2}|^{2}$ are positive since $x$ is non-degenerate. Thus we can write
$|x^{1}|^{2}=a^{2}>0$, $|x^{2}|^{2}=b^{2}>0$. It then follows that $M^{2}=\mathbb{S}^{1}(a)\times\mathbb{S}^{1}(b)$.

Finally, by the assumption \eqref{1.4},
$|h|^{2}\geq 2$, it should holds that $a^{2}+b^{2}\geq 2a^{2}b^{2}$.\endproof

\end{document}